\documentclass[12pt]{amsart}

 \usepackage{amsfonts,amssymb,eucal}
\usepackage{amsthm} 
 \usepackage{amsmath} 
\usepackage{amscd}
 \usepackage{latexsym} 
\input{cyracc.def}

\numberwithin{equation}{section}
\usepackage{color}
 \def\ii{\operatorname{i}}
\newcommand{\dd}{\operatorname{d}}
\newcommand{\Ad}{\ensuremath{{\mbox{\rm{Ad}}}}}

\newcommand{\e}{{\mbox{\rm e}}}

\newcommand{\mb}[1]{{\mbox{\boldmath{$#1$}}}}
\newcommand{\mc}[1]{{\mathcal{#1}}}
\newcommand{\got}[1]{{\mathfrak{#1}}}
\newcommand{\db}[1]{{\mathbb{#1}}}

\newcommand{\pa}{\partial}
\newcommand{\R}{\ensuremath{\mathbb{R}}}
\newcommand{\C}{\ensuremath{\mathbb{C}}}

\newcommand{\Hi}{\ensuremath{\got{H}}}

\newtheorem{deff}{Definition}
\newcommand{\h}{\ensuremath{\got{h}}}
\newcommand{\m}{\ensuremath{\got{m}}}

\newcommand{\g}{\ensuremath{\got{g}}}

\newcommand{\Z}{\ensuremath{\mathbb{Z}}}
\newcommand{\Hinf}{\ensuremath{\mathcal{H}^{\infty}}}

\renewcommand{\P}{\ensuremath{\mathbb{P}}}
\newcommand{\Phinf}{\ensuremath{\P (\Hinf )}}
\newtheorem{Remark}{Remark}

\newcommand{\D}{\ensuremath{{\got{D}}}}
\newcommand{\fl}{\ensuremath{{\got{F}}_{\Hi}}}

\newtheorem{Proposition}{Proposition}
\newtheorem{lemma}{Lemma}

\theoremstyle{definition}

\begin{document}
\title{Coherent states  and geometry on the Siegel-Jacobi disk} 
\author{Stefan  Berceanu}
\address[Stefan  Berceanu]{National
 Institute for Physics and Nuclear Engineering\\
         Department of Theoretical Physics\\
         PO BOX MG-6, Bucharest-Magurele, Romania}
\email{Berceanu@theory.nipne.ro}

\begin{abstract}
The coherent state  representation of the Jacobi group $G^J_1$ is 
indexed  with two parameters, $\mu (=\frac{1}{\hbar})$,  describing 
the part coming from the Heisenberg group, and $k$,  characterizing 
 the positive discrete series representation of $\text{SU}(1,1)$.  The
 Ricci form, the scalar curvature and the  geodesics of the Siegel-Jacobi
 disk $\mc{D}^J_1$ are investigated. The significance in the language
 of coherent states of  the transform which realizes the fundamental conjecture
 on the Siegel-Jacobi disk is emphasized. The Berezin kernel, Calabi's
 diastasis, the Kobayashi embedding, and the Cauchy formula for the
 Sigel-Jacobi disk are presented. 
\end{abstract}
\maketitle
\noindent
\tableofcontents
\newpage
\section{Introduction}
\enlargethispage{1cm}
The Jacobi group of index $n$,   $G^J_n$,   is the semidirect product of the 
symplectic group $\text{Sp}(n,\R)$ with the real $(2n+1)$-dimensional
Heisenberg group  $H_n$ \cite{ez,bs}. The points of the homogeneous
domain $\mc{D}^J_n$  associated to the Jacobi group
$G^J_n$,
called the Siegel-Jacobi domain \cite{Y07,Y08},  are in $\mc{D}_n\times\C^n$, where
$\mc{D}_n$ denotes the Siegel ball. The Jacobi groups are unimodular,
nonreductive, algebraic  groups, and the
Siegel-Jacobi spaces are  reductive,
nonsymmetric domains associated to the Jacobi groups \cite{Y07,Y08,gem} by the
generalized Harish-Chandra embedding \cite{satake}.  
The holomorphic irreducible unitary representations of the Jacobi groups
based on Siegel-Jacobi domains have been constructed \cite{bb,bs,TA90,TA92,TA99,gem}.

Among  many other applications \cite{jac1,sbj,sbcg,gem,nou,FC},  the Jacobi
group is  responsible  for the  {\it squeezed states} \cite{ken,stol,lu,yu,ho} in quantum optics
\cite{mandel,ali,siv,dr}.

We have associated to the Jacobi group $G^J_n$ coherent states in the
meaning of Perelomov \cite{perG}, based on the Siegel-Jacobi space
$\mc{D}^J_n$
\cite{sbj,nou}. Similar constructions have been considered in
\cite{KRSAR82, Q90}.  In \cite{jac1} we have defined coherent states
based on the Siegel-Jacobi disk $\mc{D}^J_1$. In the present paper we revisit  our previous
construction in \cite{jac1}, introducing a parameter $\mu\in\R$ to
describe  the holomorphic representation associated with the Heisenberg
group $H_1$, while the integer $k$ indexes  
 the
holomorphic positive discrete series of $\text{SU}(1,1)$
\cite{bar47}.
The standard realization of the position-momentum operators ${\hat{\mb{q}}}=q , ~~{\hat{\mb{p}}}=-\ii\hbar\frac{\pa}{\pa q} $ in quantum
mechanics corresponds to $\mu=\frac{1}{\hbar}$ \cite{MD}.

 In this paper we look at several
results obtained in our previous publications devoted to the Jacobi
group \cite{jac1}-\cite{FC} from the point of view of a program of
 investigation of differential geometry using the coherent states
defined   on symmetric domains \cite{ber94}-\cite{bertr},  this time considering coherent states
based on the Siegel-Jacobi spaces.

In order to make the paper selfcontained, the notation  adopted for coherent states is briefly recalled in
\S \ref{s22}. More details are given in \cite{sb6,jac1}. Our point of view on the deep relationship between coherent
states and geometry is summarized in \S \ref{geom}.  In \S \ref{jac1}
we 
collect  basic facts on  the Jacobi algebra $\got{g}^J_1$ and the Jacobi
group $G^J_1$. Using the
Bargmann transform,  Remark
\ref{rem11} determines  $\mu=\frac{1}{\hbar}$ in order to have a
representation of the Heisenberg group  compatible with the standard quantum mechanics. The $FC$-transform, 
 relating the normalized and the unnormalized Perelomov coherent
state vectors based on the Siegel-Jacobi disk $\mc{D}^J_1$, which has
an important significance below, is written down explicitly in Lemma \ref{lema6}. Most of
the result in \S \ref{OFG}, devoted to the geometry of the
Siegel-Jacobi disk, are new. From the homogeneous  metric $\dd
s^2(z,w)_{k\mu}$, 
we determine in Proposition \ref{prop2} the Ricci form, which depends only on
$w\in\mc{D}_1$, and  the Siegel-Jacobi disk is not an Einstein manifold.  The scalar
curvature of $\mc{D}^J_1$ is constant and negative. Similar results where obtained
also in \cite{jae}, working on the Sigel-Jacobi upper half plane. In
Proposition \ref{prop3}
it is observed that the change of coordinates that realizes the
fundamental conjecture
 \cite{GV,gpv,DN} for the Siegel-Jacobi disk \cite{FC1,FC,nou} has a 
significance in the context of coherent states. A discussion of the 
geodesics of the Siegel-Jacobi disk  is presented in \S \ref{GEOS}. The Berezin kernel, Calabi's diastasis, the Kobayashi embedding, and the Cauchy formula for the Sigel-Jacobi disk are
calculated  in  \S\ref{EMBSS}. The paper ends with
an appendix containing the definition of the reductive homogeneous
spaces and a few other notions of differential geometry which appear
in the body of the paper. In \cite{FC1,nou, FC} we have presented some physical
systems based on the Siegel-Jacobi spaces described by coherent states
attached to the Siegel-Jacobi group.

In this paper the Hilbert space $\Hi$ is endowed with a scalar product
antilinear in the first argument, $<\lambda x,y>=\bar{\lambda}<x,y>$,
$x,y\in\Hi,\lambda\in\C\setminus 0$. We denote the imaginary unit
$\sqrt{-1}$ by $\ii$, and the Real and Imaginary part of a complex
number by $\Re$ and respectively $\Im$, i.e. we have for $z\in\C$,
$z=\Re z+\ii \Im z$. 

\section{Coherent states - an introduction}\label{s22}
In order to fix the  notation on coherent states \cite{perG}, let us consider the triplet $(G, \pi , \got{H} )$, where $\pi$ is
 a continuous, unitary, irreducible 
representation
 of the  Lie group $G$
 on the   separable  complex  Hilbert space \Hi . 

For $X\in\g$, where \g ~is the Lie algebra of the Lie group $G$,  let
us define the (unbounded) operator $\dd\pi(X)$ on \Hi~ by
$\dd\pi(X).v:=\left. {\dd}/{\dd t}\right|_{t=0} \pi(\exp tX).v$,
whenever the limit on the right hand side exists. 
   We obtain a
representation of the Lie algebra \g~ on the smooth vectors   \Hinf of
\Hi,  {\it the derived
representation}, and we denote
${\mb{X}}.v:=\dd\pi(X).v$ for $X\in\g ,v\in \Hinf$. 

Let us now denote by $H$  the isotropy group with Lie subalgebra
$\got{h}$ of the Lie algebra $\got{g}$.  
We  consider (generalized) coherent
 states on complex  homogeneous manifolds $M\cong
G/H$ \cite{perG}.  We consider manifolds $M$  which are CS-{\em
  orbits}, i.e.  which admit a holomorphic
embedding   $\iota_M : M \hookrightarrow \Phinf$ \cite{lis,neeb,sb6}.  



 It can be shown \cite{last} that  {\em the
CS-manifolds  are  reductive spaces}. Let us denote by $\got{m}$ the
vector space orthogonal to the Lie algebra $\got{h}$ of $H$ in
$\got{g}$, i.e. we have the vector space decomposition
$\got{g}=\got{h}+\got{m}$. 
{\it The tangent space to $M$ at $o$
can be identified with $\got{m}$} and if $\exp : \got{g}\rightarrow G$
is the exponential mapping, then $G/H=\exp(\got{m})$,  see definition
at p. 104 in \cite{helg} and the Appendix. 

We
can introduce the normalized (unnormalized) vectors  $\underline{e}_x$
(respectively, $e_z$) defined on $G/H$
\begin{equation}\label{bigch}
\underline{e}_x=\exp(\sum_{\phi\in\Delta_+}x_{\phi}\mb{X}_{\phi}^+-\bar{x}_{\phi}\mb{X}_{\phi}^-)
e_0, ~
e_z=\exp(\sum_{\phi\in\Delta_+}z_{\phi}\mb{X}_{\phi}^+)e_0, 
\end{equation}
where
$e_0$ is the extremal weight vector of the representation $\pi$, $\Delta_+$ are the positive roots of the Lie algebra $\got{g}$ of $G$,
and $X_\phi,\phi\in\Delta$,  are the  generators. 
$\mb{X}^+_{\phi}$ ($\mb{X}^-_{\phi}$)
corresponds to the positive (respectively, negative) generators. See
details in  \cite{perG,sb6}.

Let us denote by $FC$ the change of variables
$x\rightarrow z$ in formula \eqref{bigch} such that
\begin{equation}\label{etild}\underline{e}_{x}=\tilde{e}_z, ~    \tilde{e}_z  :=
  (e_z,e_z)^{-\frac{1}{2}}e_z, ~z=FC(x). 
\end{equation}
The reason for calling the transform \eqref{etild}  $FC$  {\it
  (fundamental conjecture)} is explained later, see Proposition \ref{prop3}.
 
{\em The coherent vector
 mapping}  $\varphi$ is defined locally, on a coordinate neighborhood
$\mc{V}_0$ (cf. \cite{last,sb6}):
 \begin{equation}\label{ECFI}\varphi : M\rightarrow \bar{\Hi}, ~
 \varphi(z)=e_{\bar{z}},
\end{equation}
where $ \bar{\Hi}$ denotes the Hilbert space conjugate to $\Hi$.
The  vectors $e_{\bar{z}}\in\bar{\Hi}$ indexed by the points
 $z \in M $ are called  {\it
Perelomov's coherent state vectors}. Using Perelomov's coherent
vectors, we consider Berezin's approach to quantization of K\"ahler
manifolds with the supercomplete sets of vectors
\cite{ber73}-\cite{berezin} in the formulation of Rawnsley \cite{raw,Cah,cah}.

Let us denote by $\fl= L^2 _{\text{hol}} (M,\dd \nu_M)$  the space of holomorphic, square integrable  functions 
  with respect to  the scalar product on $M$
\begin{equation}\label{scf}
(f,g)_{\fl} =\int_{M}\bar{f}(z)g(z)\dd \nu_M(z,\bar{z}),
~~\dd{\nu}_{M}(z,\bar{z})=\frac{\Omega_M(z,\bar{z})}{K_M}, ~K_M= (e_{\bar{z}},e_{\bar{z}}).
\end{equation}

The positive real function $K_M=K_M(z,\bar{z})$  in \eqref{scf}, 
also  denoted $K(z)$, is the Bergman kernel, expressed as the scalar
product of coherent states based on the K\"ahler  homogenous manifold $M=G/H$, 
 $\Omega_M$ is the normalized  $G$-invariant volume form
\begin{equation}
\Omega_M:=(-1)^{\binom {n}{2}}\frac {1}{n!}\;
\underbrace{\omega\wedge\ldots\wedge\omega}_{\text{$n$ times}}\ ,
\end{equation}
and the $G$-invariant  K\"ahler two-form $\omega$ on the $2n$-dimensional
manifold  $M=G/H$  with K\"ahler potential 
\begin{equation}\label{FK}
f(z,\bar{z})=\ln K(z,\bar{z})
\end{equation} 
is given by
\begin{equation}\label{kall}
\omega_M(z)=\ii\sum_{\alpha,\beta\in\Delta_+} h_{\alpha\bar{\beta}} (z) \dd z_{\alpha}\wedge
\dd\bar{z}_{\beta}, 
\end{equation}
\begin{equation}\label{Ter} 
h_{\alpha\bar{\beta}}= \frac{\pa^2 f}{\pa {z}_{\alpha}\pa \bar{z}_{{\beta}}} = \frac{\pa^2}{\pa
  z_{\alpha} \pa\bar{z}_{\beta}} \ln (e_{\bar{z}},e_{\bar{z}}),
~h_{\alpha,\bar{\beta}}= \bar{h}_{\beta\bar{\alpha}}.
\end{equation}
 The {\it hermitian  (Bergman) metric} of $M$ in
local coordinates  is 
\begin{equation}\label{herm}
\dd s^2_M(z,\bar{z})=\sum_{\alpha,{\beta}}h_{\alpha\bar{\beta}}\dd
z_{\alpha}\otimes \dd\bar{z}_{\beta} =\sum_{\alpha,{\beta}} \frac{\pa^2}{\pa
  z_{\alpha} \pa\bar{z}_{\beta}} \ln (K_M(z,\bar{z}))
\dd
z_{\alpha}\otimes \dd\bar{z}_{\beta} .
\end{equation} 
and the
condition of the metric to be a K\"ahlerian one is (cf. (6) p. 156 in \cite{koba2})
\begin{equation}\label{condH}
\frac{\pa h_{\alpha\bar{\beta}}}{\pa z_{\gamma}} = \frac{\pa
  h_{\gamma\bar{\beta}}}{\pa z_{\alpha}},~\alpha,\beta,\gamma
=1,\dots, n .
\end{equation}
If $\left\{\varphi_n(z)\right\}_{n=1,\dots}$ is an orthonormal  base of functions
of \fl  and \begin{equation}\label{funck}
K_M(z,\bar{w}):=(e_{\bar{z}},e_{\bar{w}}), \end{equation} then the
Bergman kernel admits the series expansion 
\begin{equation}\label{sumBERG}
K_M(z,\bar{w})=\sum_n^{\infty}\varphi_n(z)\bar{\varphi}_n(w).
\end{equation}
For compact manifolds $M$, \fl~ is finite dimensional  and also
the sum \eqref{sumBERG} is finite.

Let us  introduce the map $\Phi :\Hi^{\star}\!\rightarrow \fl$ ,
\begin{equation}\label{aa}
\Phi(\psi):=f_{\psi},
f_{\psi}(z)=\Phi(\psi )(z)=(\varphi (z),\psi)_{\Hi}=(e_{\bar{z}},\psi)_{\Hi},
~z\in{\mathcal{V}}_0\subset M, 
\end{equation}
where we have identified the space $\overline{\Hi}$  complex conjugate
 to \Hi~  with the dual
space
$\Hi^{\star}$ of $\Hi$.  \eqref{scf}
is nothing else but   {\em {Parseval   overcompletness
 identity}}  \cite{berezin}
\begin{equation}\label{overc}
(\psi_1,\psi_2)=\int_{M=G/H}(\psi_1,e_{\bar{z}})(e_{\bar{z}},\psi_2)\dd
\nu_M(z,\bar{z}). 
\end{equation}

The relation \eqref{scf} can be interpreted in the language of
geometric quantization \cite{Kos}. Together with the K\"ahler manifold
$(M,\omega)$, we also consider the triple $\sigma =(L,h,\nabla)$, where
$L$ is a holomorphic line bundle on $M$, $h$ is the hermitian metric
on $L$ and $\nabla$ is a connection compatible with metric and the
K\"ahler structure \cite{SBS}. With respect to a local holomorphic
frame for the line  bundle, the metric can be given as
$$h(s_1,s_2)(z)=\hat{h}(z)\bar{\hat{s}}_1(z)\hat{s}_2(z),$$
where $\hat{s}_i$ is a local representing function for the section
$s_i$, $i=1, 2$,  and $\hat{h}(z)=(e_z,e_z)^{-1}$. \eqref{scf} is the local
representation of the scalar product on the line bundle
$\sigma$. The connection $\nabla$ has the expression  $\nabla=\pa +\pa \ln \hat{h} +\bar{\pa}$. The curvature of
$L$ is defined as
$F(X,Y)=\nabla_X\nabla_Y-\nabla_Y\nabla_X-\nabla_{[X,Y]}$, and locally
$F=\bar{\pa}\pa\ln\hat{h}$ \cite{chern}. The K\"ahler manifold  $(M,\omega)$  is
{\it quantizable} if there exists a triple $\sigma$ such that
$F(X,Y)=-\ii\omega(X,Y)$ and  we have \eqref{Ter}.

The square of the length of a vector $X\in\C^n$, measured in this
metric at the point $z\in M$,  is \cite{dinew}
$$\tau_M^2(z,X):=\sum_{\alpha,\beta}^nh_{\alpha\bar{\beta}}(z)X_{\alpha}\bar{X}_{\beta} .$$
If the length $l$ of a piecewise $C^1$-curve $\gamma:[0,1]\ni t\mapsto
\gamma(t)\in M$ is defined as 
$$l(\gamma):=\int_0^1\tau_M(\gamma(t),\gamma'(t)) \dd t,$$
then the {\it Bergman distance} between two points $z_1,z_2\in M$ is
$$\dd_B(z_1,z_2)=\text{inf}\left\{ l(\gamma):\gamma ~\text{is a
      piecewise curve s.t.} ~\gamma(0)=z_1,\gamma(1)=z_2\right\}.$$

We  denote  the ``normalized'' Bergman kernel by 
\begin{equation}\label{kmic}
\kappa_M(z,\bar{z}'):=\frac{K_M(z,\bar{z}')}{\sqrt{K_M(z)K_M(z')}}=
(\tilde{e}_{\bar{z}},\tilde{e}_{\bar{z}'})=\frac{(e_{\bar{z}},e_{\bar{z}'})}{||e_{\bar{z}}||||e_{\bar{z}'}||} .
\end{equation}
Introducing in the above definition  the series expansion
\eqref{sumBERG}, with  the Cauchy - Schwartz inequality, we have that 
$$|\kappa_M(z,\bar{z}')|\le 1.$$
In this paper  by  {\it  the Berezin kernel} $b_M:M\times M\rightarrow
[0,1]\in \R$ we mean:
\begin{equation}\label{berK}
b_M(z,z')=|\kappa_M(z,\bar{z}')|^2.
\end{equation}

\section{The coherent states - a bridge between quantum mechanics and geometry}\label{geomp}
 
In  \cite{ber94} we have advanced the proposal of a program  of
investigation of the  deep relationship  between  coherent
states and the geometry of the manifolds on which the coherent
states are defined. We have investigated how the
coherent states permit to find:  1)
the geodesics; 2) the conjugate locus; 3) the cut locus; 4) the
divisors; 5) the Calabi's diastasis. Also, we wanted to find a
geometric meaning of transition amplitudes  and of different distances and angles in quantum mechanics for
coherent states.   The obtained  results are
established  for
very particular manifolds, mostly on  hermitian symmetric spaces \cite{ber94}-\cite{bertr}. Our favorite
example was  the complex Grassmann manifold and its noncompact dual
\cite{gr}.

In \cite{ber94,ber95Q,cut,bertr} we have proved that
\begin{Remark}\label{trecut} 
For symmetric spaces
 the dependence $z(t)=FC(tx)$  from \eqref{etild} gives geodesics in $M$ with the
property that $z(0)=p$ and $\dot{z}(0)=x$.  
\end{Remark}

Let $\xi:\Hi\setminus 0\rightarrow\db{P}(\Hi)$ be the  the  canonical projection 
$\xi(\mb{z})=[\mb{z}]$. The {\it Fubini-Study metric}  in the
nonhomogeneous coordinates $[z]$ is  the
hermitian metric on $\db{CP}^{\infty}$ \cite{koba}
\begin{equation}\label{FBST}\dd s^2|_{FS}(\mb[{z}])=
\frac{(\dd\mb{z},\dd\mb{z}) (\mb{z},\mb{z})-(\dd\mb{z},\mb{z})
  (\mb{z},\dd\mb{z})}{(\mb{z},\mb{z})^2}.
\end{equation}

The elliptic {\it Cayley distance}    \cite{Cay}  between two points in the projective Hilbert space
$\db{P}(\Hi)$ is defined as 
\begin{equation}\label{Cdis}
\dd_C([z_1],[z_2])=\arccos \frac{|(z_1,z_2)|}{||z_1|| ||z_2 || }.
\end{equation}
The  Fubini-Study metric \eqref{FBST} and the Cayley distance
\eqref{Cdis} are  independent of the homogeneous coordinates $z$
representing $[z]=\xi(z)$. 

{\it Calabi's diastasis} \cite{calabi}, in the context of coherent states
 as used by Cahen, Gutt and Rawnsley \cite{cah}, reads: 
\begin{equation}\label{DIA}D_M(z,z')= -\ln b_M(z,z') = -2\ln
\left\vert(\tilde{e}_{\bar{z}},\tilde{e}_{\bar{z}'})\right\vert.
\end{equation}

Let $M$ be a homogeneous K\"ahler manifold $M=G/H$ to which we
associate the Hilbert space of functions \fl~ with respect to the scalar
product \eqref{scf}. 

We can make the following assertions
\begin{Remark}\label{HTR}Let us suppose that the K\"ahler manifold $M$ admits a holomorphic
embedding 
\begin{equation}\label{invers}
\iota_M: M\hookrightarrow \db{CP}^{\infty}, \iota_M(z) =  [\varphi_0(z):\varphi_1(z):\dots ] .
\end{equation}
The Hermitian metric \eqref{herm} on $M$ is
the pullback of the Fubini-Study metric \eqref{FBST} via the embedding
\eqref{invers}, i.e.:
\begin{equation}\label{KOL} 
\dd s^2_M(z)=\iota_M^*\dd s^2_{FS}(z)= \dd s^2_{FS}(\iota_M(z)).
\end{equation}  
The angle defined by the normalized Bergman kernel \eqref{kmic} can
be expressed via the embedding \eqref{invers}  as function of the Cayley
distance \eqref{Cdis}
\begin{equation}\label{c1}
\theta_M(z_1,z_2)=\arccos|\kappa_M(z_1,\bar{z}_2)|=
\arccos|(\tilde{e}_{z_1},\tilde{e}_{z_2})_M|
=d_C(\iota_M(z_1),\iota_M(z_2)).
\end{equation}

We have also the relation 
\begin{equation}\label{ddbm}\dd_B(z_1,z_2)\ge\theta_M(z_1,z_2) . 
\end{equation}
The following (Cauchy) formula is true
\begin{equation}\label{c2}
 (\tilde{e}_{z_1},\tilde{e}_{z_2})_M=(\iota_M(z_1),\iota_M(z_2))_{\db{CP}^{\infty}}.
\end{equation}
The Berezin kernel \eqref{berK} admits the geometric interpretation via the Cayley
distance as
\begin{equation}
b_M(z_1,z_2)=\cos^2d_C(\iota_M(z_1),\iota_M(z_2))=\frac{1+\cos(2d_C(\iota_M(z_1),\iota_M(z_2)))}{2} .
\end{equation}
\end{Remark}
\begin{proof}
The assertion \eqref{KOL} is known from \cite{koba}. 
We introduce in the expression \eqref{FBST} of the Fubini-Study
hermitian metric
on $\db{CP}^{\infty}$  the change of coordinates realized by the
embedding $\iota_M$ \eqref{invers} and we get
\begin{equation}\label{sumJK} \dd s_{FS}^2(z)=\frac{(\sum\dd \bar{\varphi}_i\dd
  \varphi_i)(\sum\bar{\varphi}_i\varphi_i)-
|\sum
\dd\bar{\varphi_i}\varphi_i|^2}{(\sum\bar{\varphi}_i\varphi_i)^2}.
\end{equation}
But the Bergman kernel $K$ admits  the
expansion \eqref{sumBERG}
$$K(z,\bar{z})=\sum\bar{\varphi}_i(z)\varphi_i(z),$$
and  we write down \eqref{sumJK} in the coordinates $z$ of the
manifold $M$ as 
\begin{equation}\label{FXKK} 
\dd s_{FS}^2(z)=\frac{K\dd_z\dd_{\bar{z}}K - (\dd_z
  K)(\dd_{\bar{z}}K)}{K^2}.
\end{equation}
In the above expression we take into account \eqref{FK},
i.e. $K=\e^f$, which implies
\begin{equation}\label{SDF}
\begin{split}
\dd_{\bar{z}}K & = K\sum\frac{\pa f}{\pa \bar{z}_j}\dd\bar{z}_j,\\
\dd_z\dd_{\bar{z}}K & =K\sum\frac{\pa f}{\pa z_k}\frac{\pa f}{\pa
  \bar{z}_j}\dd z_k\dd\bar{z}_j+K\sum\frac{\pa ^2f}{\pa z_k\pa
  \bar{z}_j}\dd z_k \dd \bar{z}_j .
\end{split}
\end{equation}
Now we introduce the expressions \eqref{SDF} into \eqref{FXKK} and we
get
$$\dd s^2_{FS}(z)=\sum_{i,j}\frac{\pa^2 f}{\pa z_i \pa \bar{z}_j } \dd z_i\dd \bar{z}_j,$$
i.e. \eqref{KOL}. The formulas \eqref{c1}, \eqref{c2} are proved in
\cite{ber94}, \cite{ber95Q}, \cite{bertr}. 
The inequality \eqref{ddbm} is proved in \cite{Block}, \cite{LU09}.
\end{proof}

\section{Coherent states on the Siegel-Jacobi disk}\label{jac1}
\subsection{The generators of the Lie algebra $\got{g}^J_1$}

The Heisenberg   group  is the group with the
3-dimensional real  Lie algebra isomorphic to the Heisenberg  algebra
 \begin{equation}\label{nr0}\got{h}_1\equiv
<\ii s 1+\alpha a^{\dagger}-\bar{\alpha}a>_{s\in\R ,\alpha\in\C} ,\end{equation}
 where $a^{\dagger}$ ($a$) are  the boson creation
(respectively, annihilation)
operators which verify the canonical commutation relations 
(\ref{baza1}).

Let us also consider the Lie algebra of the group $\text{SU}(1,1)$:
\begin{equation}\label{nr1}
\got{su}(1,1)=
<2\ii\theta K_0+yK_+-\bar{y}K_->_{\theta\in\R ,y\in\C} , \end{equation} 
where the generators $K_{0,+,-}$ verify the standard commutation relations
(\ref{baza2}).

 Now let us define the Jacobi algebra as the  the semi-direct sum \cite{jac1}
\begin{equation}\label{baza}
\got{g}^J_1:= \got{h}_1\rtimes \got{su}(1,1),
\end{equation}
where $\got{h}_1$ is an  ideal in $\got{g}^J_1$,
i.e. $[\got{h}_1,\got{g}^J_1]=\got{h}_1$, 
determined by the commutation relations \eqref{baza3}, \eqref{baza5}:
\begin{subequations}\label{baza11}
\begin{eqnarray}
& & [a,{{a}^\dagger}]=1\label{baza1}, \\
\label{baza2}
~& & \left[ K_0, K_{\pm}\right]=\pm K_{\pm}~,~ 
\left[ K_-,K_+ \right]=2K_0 , \\
\label{baza3}
& & \left[a,K_+\right]=a^{\dagger}~,~\left[ K_-,a^{\dagger}\right]=a, ~
\left[ K_+,a^{\dagger}\right]=\left[ K_-, a\right]= 0 ,\\
\label{baza5}
& & \left[ K_0  ,~a^{\dagger}\right]=\frac{1}{2}a^{\dagger}, \left[ K_0,a\right]
=-\frac{1}{2}a .
\end{eqnarray}
\end{subequations}
\begin{Remark}\label{geom}
The Jacobi groups are unimodular, non-reductive, algebraic groups of
Harish - Chandra type. The Siegel-Jacobi domains are reductive,
non-symmetric manifolds associated to the Jacobi groups by the
generalized Harish-Chandra embedding.
\end{Remark}
\begin{proof}
From the algebraic relations \eqref{baza11}, we can conclude that
{\it the Sigel-Jacobi disk  $\mc{D}^J_1$
is a non-symmetric space}. Indeed,  $\mc{D}^J_1=G^J_1/H$, where 
$H=U(1)\times\C$ is connected  and the Jacobi algebra admits the
decomposition of the type \eqref{dec}, i.e. $\got{g}^J_1=\got{h}+\got{m}$, where
$\got{h}$ is generated by $K_0$ and $1$. We have
$[\got{h},\got{h}]\subset\got{h}$, $[\got{h},\got{m}]\subset\got{m}$,
and {\it the Siegel-Jacobi domains are reductive homogeneous spaces}, in
sense of  Nomizu \cite{nomizu}, but because $[\got{m},\got{m}]\subset\got{h}$ is not
true,   the Siegel-Jacobi disk  is not a 
symmetric  space. Also, because of the presence of the Heisenberg
subalgebra in the Jacobi algebra, {\it the Jacobi group is a
  non-reductive one}. 

From  the definition of the Jacobi algebra $\got{g}^J_1$, we see
that {\it   the Jacobi group is an algebraic group of Harish-Chandra 
  type} \cite{satake}, as was underlined in \cite{gem} for $\g^J_n$. Indeed, in the
direct sum of vector spaces,
$\got{g}^J_1=\got{p}_+ +\got{k}_{\C}+\got{p}_-$, with properties 
$  [\got{k}_{\C},\got{p}_{\pm}]\subset \got{p}_{\pm}$, and
$\bar{\got{p}}_+=\got{p}_-$, we have the identification
$\got{p}_+\equiv <K_+, a^{\dagger}>$, $\got{p}_-\equiv <K_, a>$, and
$\got{k}_{\C} \equiv  <K_0,1 > $,  see also  \cite{jac1}.  The generalized Harish-Chandra embedding of
the homogeneous space $\mc{D}^J_1$ into $\got{p}_+$   is determined by
$gK\mapsto \zeta = (z,w)$, where $\exp \zeta=(g)_+$, as explained in
detail \cite{gem}, and the $G^J_1$-invariant structure of $\mc{D}^J_1$
is determined by the natural inclusion $\mc{D}^J_1\hookrightarrow
P_+=\exp\got{p}_+$. 
\end{proof}


\subsection{Formulas for the Heisenberg group $H_1$}\label{unul}
The displacement operator  $D_{\mu}(\alpha
)$ (or simply, $D(\alpha )$), $\mu\in\R_+$:
\begin{equation}\label{deplasare}
D_{\mu}(\alpha )=\exp \sqrt{\mu}(\alpha {\mb{a}}^{\dagger}-\bar{\alpha}{\mb{a}})=\exp(-\frac{1}{2}\mu|\alpha
|^2)  \exp (\sqrt{\mu} \alpha {\mb{a}}^{\dagger})\exp(-\sqrt{\mu}\bar{\alpha}{\mb{a}}),
\end{equation}
has the composition property
\begin{equation}\label{thetah}
D_{\mu}(\alpha_2)D_{\mu}(\alpha_1)=e^{\ii\theta_{\mu}(\alpha_2,\alpha_1)}
D_{\mu}(\alpha_2+\alpha_1) , 
~\theta_{\mu}(\alpha_2,\alpha_1):=\mu\Im (\alpha_2\bar{\alpha_1}) .
\end{equation}
The coherent states associated to the Heisenberg group $H_1$ are
defined on the homogeneous space $\C=H_1/\R$ taking in \eqref{etild}
as $e_0$ the vacuum $e^H_0$, with $\mb{a}e^H_0=0$ and then
\eqref{bigch} has the form
$$ \underline{e}_{\alpha}:=  D_{\mu}(\alpha)e_0=   \e^{-\frac{1}{2}\mu|\alpha|^2} e_{\alpha},
 ~ e_{\alpha}:=\e^{\sqrt{\mu}\alpha\mb{a}^{\dagger}}e_0, $$
a particular case of the relation
\begin{equation}\label{rePP}   D_{\mu}(\alpha)e_z=
  \e^{-\mu\bar{\alpha}(z+\frac{\alpha}{2})} e_{\alpha+z}.
\end{equation}
Equation \eqref{rePP} corresponds to the projective representation of the parallel
transport in the Hilbert space $\got{F}_{\mu}$  with scalar product
\eqref{scfFF}; see  equation (4.17) in \cite{ber74}, or equation (107)
in \cite{cartier66};  see p. 47 in 
\cite{fol}. For the   Fock representation 
$$\beta_{\mu}(\alpha,t):=\e^{\ii\mu t}D_{\mu}(\alpha)$$ of the
Heisenberg algebra  \eqref{nr0}  we have to take $\dd \beta_{\mu} (1) =\ii\mu$,
and 
$$\beta_{\mu}(\alpha,t)e_z=\e^{\ii\mu
  t}\e^{-\mu\bar{\alpha}(z+\frac{\alpha}{2})} e_{\alpha+z},$$
where $\mu=2\pi m,~ m\in\R$. In \cite{jac1} we have taken $\mu=1$.

The orthonormal base of $n$-vectors consists of the vectors
\begin{equation}\label{ort2}
|n>:=(n!)^{-\frac{1}{2}}({\mb{a}}^{\dagger})^{n}|0>;~~<n',n>=\delta_{nn'} .
\end{equation} 
  Perelomov's CS-vectors associated to the Heisenberg
group, defined on $M:=H_1/\R=\C$, are
\begin{equation}\label{v1}
 e_z:=\e^{\sqrt{\mu}z{\mb{a}}^{\dagger}}e_0=\sum\frac{(\sqrt{\mu}z)^n}{(n!)^{1/2}}|n> .
\end{equation}
The  corresponding holomorphic functions  associated to \eqref{v1}  are (see e.g. formula (4.3) in \cite{ber74}, or Theorem 3.4 in \cite{hall00}, or  equation (112) in \cite{cartier66}; Bargmann took $\mu=1$  in equation (1.6)  in \cite{bar} as we took in our previous publication \cite{jac1}):
\begin{equation}\label{f1}
f_{|n>}(z):=(e_{\bar{z}},|n>)=\frac{(\sqrt{\mu}z)^n}{(n!)^{1/2}} .
\end{equation} 

The reproducing kernel $K_{\mu}:\C\times\bar{\C}\rightarrow\C$ is
\begin{equation}\label{ker1}
K _{\mu} (z,\bar{z}'):=(e_{\bar{z}},{e}_{\bar{z}'})=\sum f_{|n>}(z)\bar{f}_{|n>}(z')=
e^{\mu z\bar{z}'} ,
\end{equation}
where the vector $e_z$ is given by  (\ref{v1}), while the function
$f_{|n>}(z)$ is given by  (\ref{f1}). 

The homogeneous space $\C=H_1/R$ has  the $H_1$-invariant  K\"ahler two-form
$\omega_{\mu}$
$$-\ii \omega_{\mu}(z)=\mu \dd z \wedge \dd \bar{z}.$$

The scalar product (\ref{scf})
on the Segal-Bargmann-Fock space  
$\got{F}_{\mu}=L^2 _{\text{hol}}(\C,\rho_{\mu})$  (in  \cite{bar}
Bargmann works with $\mu =1$) is 
\begin{equation}\label{scfFF}
(f,g)_{\mu}=\int_{\C}\bar{f}g\rho_{\mu} \dd\Re z \dd\Im z , ~~~\rho_{\mu}=\frac{\mu}{\pi}\exp{(-\mu|z|^2)}.
\end{equation}

Note that 
\begin{Remark}\label{rem11}
 Under  the differential realization on $\got{F}_{\mu}$ of the creation and annihilation operators
\begin{equation}\label{eqr}
\mb{a}^{\dagger}=\sqrt{\mu}z, ~~\mb{a}=\frac{1}{\sqrt{\mu}}\frac{\pa}{\pa z},
\end{equation} 
the operators $\mb{a}$ and $\mb{a}^{\dagger}$  are hermitian conjugate
with respect to  the scalar product \eqref{scfFF}. 
The standard realization
\begin{equation}\label{std}
\hat{\mb{q}}=q , ~~\hat{\mb{p}}=-\ii\hbar\frac{\pa}{\pa q} 
\end{equation} 
 in $\got{H}=L^2(\R,\dd x)$ of the position and momentum operators, 
where
\begin{equation}\label{aapl}\mb{a}=\lambda(\hat{\mb{q}}+\ii \hat{\mb{p}}), ~~
  \mb{a}^{\dagger}=  \lambda(\hat{\mb{q}}-\ii \hat{\mb{p}}) ,\end{equation} 
corresponds to the choice of $\mu$ in
\eqref{ker1}, \eqref{scfFF} and $\lambda$  in \eqref{aapl} as 
\begin{equation}\label{condS}
\mu \hbar =1, ~ ~2\hbar\lambda^2 =1.
\end{equation}
\end{Remark}
\begin{proof}
It easy to verify that in  the realization \eqref{eqr}, we have 
$$(\mb{a}f,g)_{\mu}=(f,\mb{a}^{\dagger}g)_{\mu}.$$ 
  Following 
Bargmann \cite{bar}, who considered the case $\mu=1$,  we shall
find  the kernel $B$ which  determines the (Bargmann) transform
$\mc{B}:
\got{H}= L^2(\R,\dd x)\rightarrow\got{F}_{\mu}$,
\begin{equation}\label{BRRU}
f(z)=\int B(z,q)\psi(q)\dd q,  
\end{equation}
such that if $\mc{B}: \got{H}\ni\psi\mapsto f\in\got{F}_{\mu}$, then,
in accord with the representation \eqref{eqr},  $\mb{a}^{\dagger}\psi\mapsto \sqrt{\mu}zf$ and
$\mb{a}\psi\mapsto \frac{1}{\sqrt{\mu}}\frac{\pa }{\pa z} f. $ But 
\begin{subequations}\label{dublu1}
\begin{align}
 \int \mc{B}(\mb{a}^{\dagger}\psi)\dd q &  = \int(\mb{a}\mc{B})\psi\dd q =\sqrt{\mu}zf= \sqrt{\mu}\int zB\psi\dd q;\\
\int\mc{B}(\mb{a}\psi)\dd q  & = \int (\mb{a}^{\dagger}\mc{B})\psi \dd q =\frac{1}{\sqrt{\mu}}\frac{\pa }{\pa z} f= \frac{1}{\sqrt{\mu}}\int \frac{\pa B}{\pa z} \psi \dd q.
\end{align}
\end{subequations}
If in \eqref{dublu1} we introduce \eqref{aapl}, \eqref{std}, we get for the kernel $B$ of the Bargmann transform the system of differential equations 
\begin{subequations}\label{dublu2}
\begin{align}
\sqrt{\mu}zB & =\lambda (q+\hbar\frac{\pa}{\pa q})B;\\
\frac{1}{\sqrt{\mu}} \frac{\pa B}{\pa z} & = \lambda (q -\hbar\frac{\pa }{\pa q}) B .
\end{align}
\end{subequations}
Solving the system \eqref{dublu2}, we get \eqref{condS}  and the kernel $B$ of the Bargmann transform \eqref{BRRU}
is  determined as
\begin{equation}\label{solBR}
B(z,q)= B_0\e^{\frac{1}{\hbar}(\sqrt{2}qz-\frac{z^2+q^2}{2})}, ~B_0= (\pi\hbar)^{-\frac{1}{4}} .
\end{equation}
The value of $B_0$ in \eqref{solBR} was  fixed choosing  the
normalization constant such that 
$$\int B(z,q) B(\bar{w},q) \dd q=K_{\frac{1}{\hbar}}(z,\bar{w}) .$$
Note also that the normalized solution in $\got{H}$ of the equations $\mb{a}\varphi_0=0$
is $\varphi_0(q)=(2\pi\hbar)^{-\frac{1}{4}}\e ^{-\frac{q^2}{2\hbar}} $ from where we obtain the values of the
n-particle vectors  \eqref{ort2} in $\got{H}$. 
\end{proof}

\subsection{Formulas for  $\text{SU}(1,1)$}\label{unulL}
Let us denote   by $S_k$ the unitary irreducible positive discrete
series representation  $D^k_+$  of the group
$\text{SU}(1,1)$ with Casimir operator $C=K^2_0-K^2_1-K^2_2=k(k-1)$,
where $k$ is the Bargmann index for $D^+_k$ \cite{bar47}. 
Let us introduce the notation $\underline{S}_k(z)=S_k(w)$, where 
 $w\in\C,~
|w|<1$ and  $z\in\C\setminus 0$, 
are related by  \eqref{u5}. We have the relations:
\begin{subequations}
\begin{eqnarray}
\underline{S}_k(z) & := &\exp (z{\mb{K}}_+-\bar{z}{\mb{K}}_-), 
;\label{u1} \\
S_k(w) & = &  \exp (w{\mb{K}}_+)\exp (\xi
{\mb{K}}_0)\exp(-\bar{w}{\mb{K}}_-);\label{u2}\\
w & = &  \frac{z}{|z|}\tanh \,(|z|), ~~\xi =\ln (1-w\bar{w}), z\not=
0,   \label{u5}
\end{eqnarray}
\end{subequations}
and  $w=0$ for $z=0$ in \eqref{u5}.  

Here we underline that, in accord with the Remark \ref{trecut}:
\begin{Remark} \label{rem66}
{\it  The
change of variables in \eqref{u5} is of the type $FC$  in the meaning
of   {\emph{Remark \ref{trecut}}} as in
\eqref{etild}, i.e.  the dependence $w(t)=w(t z) =  \frac{z}{|z|}\tanh \,(t|z|)$ gives geodesics for the symmetric space $\mc{D}_1$ 
 with the  property that $\dot{w}(0)=z$ and $w(0)=0$}.
\end{Remark} 
\begin{proof}
Indeed, we have to verify that $w(t)=\frac{z}{|z|}\tanh \,(t|z|)$
verifies the equation of geodesics  on Siegel disk $\mc{D}_1$,
$G_2=0$, where $G_2$ is given in \eqref{geo2} below. The fact that
that  $\mc{D}_1$ is a symmetric space is evident from the structure of
the algebra  $\got{su}(1,1)$ given in \eqref{baza2} which admits a
decomposition of the  type \eqref{sum1} with property \eqref{symsp}. 
\end{proof}

In order to construct coherent states associated to $\text{SU}(1,1)$
defined on  the Siegel disk $\mc{D}_1=\text{SU}(1,1)/\text{U}(1)$, in formula
\eqref{etild} we take $e_0=e^K_o\equiv e_{kk}$ as the extreme vector of the $D^+_k$
representation such that $\mb{K}_0e_{k,k}=k e_{k,k}, k>0$. 

In the  orthonormal base  (\ref{ort1}),
\begin{equation}\label{ort1}
e_{k,k+m}:=a_{km}({\mb{K}}_+)^me_{k,k};~~a^2_{km}=\frac{\Gamma
(2k)}{m!\Gamma (m+2k)},
\end{equation} 
Perelomov's CS-vectors for $\text{SU}(1,1)$, based on the unit disk
$\mc{D}_1=\text{SU}(1,1)/\text{U}(1)$,  are
\begin{equation}\label{v2}
e_z:=e^{z{\mb{K}}_+}e_0= \sum \frac{z^n{\mb{K}}_+^n}{n!}e_0= 
\sum \frac{z^n e_{k,k+n}}{n!a_{kn}} ,\end{equation}
and  the corresponding holomorphic functions 
are (see e.g.   equation 9.14 in \cite{bar47})
\begin{equation}\label{f2}
f_{e_{k,k+n}}(z):=(e_{\bar{z}},e_{k,k+n})
=\sqrt{\frac{\Gamma (n+2k)}{n!\Gamma (2k)}}z^n .
\end{equation}

The reproducing kernel $K_k:\mc{D}_1\times\bar{\mc{D}}_1\rightarrow \C$ is
\begin{equation}\label{ker2}
K _k (z,\bar{z}'):=(e_{\bar{z}},{e}_{\bar{z}'})=
\sum f_{e_{k,k+m}}(z)\bar{f}_{e_{k,k+m}}(z')=
(1-z\bar{z}')^{-2k}, 
\end{equation}
where the vector $e_z$ is given by  (\ref{v2}), while the function
$f_{e_{k,k+m}}(z)$ is given by  (\ref{f2}).

The Siegel disk $\mc{D}_1$ has the  K\"ahler two-form $\omega_k$
$$-\ii \omega_k(w) =\frac{2k}{(1-w\bar{w})^2}\dd w\wedge \dd
\bar{w},  $$
$\text{SU}(1,1)$-invariant to the linear fractional  transformation 
\begin{equation}\label{dg} g\cdot w 
=\frac{a \, w+ b}{\delta}, ~ \delta=\bar{b}w+\bar{a},
\text{SU}(1,1) \ni g= \left( \begin{array}{cc}a & b\\ \bar{b} &
\bar{a}\end{array}\right),~\text{where} ~|a|^2-|b|^2=1.
\end{equation} 

The action \eqref{dg} is {\it  transitive}, and   the Siegel disk
$\mc{D}_1$ is a  {\it
  homogeneous manifold}.  As was already mentioned in the proof of  Remark \ref{rem66}, $\mc{D}_1$  is symmetric:
the solution $g\in\text{SU}(1,1)$ of the equation  $\sigma:\mc{D}_1\ni
w\mapsto -w\in\mc{D}_1$ is $g$  in \eqref{dg} with 
$a=\pm 1, ~ b=0$. Then $\sigma^2=1$, and  $0$ is an isolated fixed
point of $\sigma$. By homogeneity,  it can be shown that for $\forall
w\in\mc{D}_1$,   there exists a $\sigma_w$ such that $\sigma^2_w=1$, and
$w$ is an isolated fixed point of $\sigma_w$.

The scalar product (\ref{scf}) on
$\mc{D}_1=\text{SU}(1,1)/\text{U}(1)$  in $\got{F}_k=L^2 _{\text{hol}}(\mc{D}_1,\rho_k)$ 
 is  (see e.g.
equation 9.9 in \cite{bar47})
\begin{equation}\label{psipsi}
(\phi ,\psi )_{k}=
\int_{|z|<1} \bar{f}_{\phi}(z)f_{\psi}(z)\rho_k\dd\Re
z\dd\Im z,~\rho_k=  \frac{2k-1}{\pi}(1-|z|^2)^{2k-2}
, ~2k= 2,3,...  .
\end{equation}
In order to prove Proposition \ref{mm1}, we recall \cite{jac1}
\begin{Remark}Let $g\in
\rm{SU}(1,1)$ has the form \eqref{dg}.

The following relations hold:
\begin{equation}\label{r1}
\underline{S}_k(z)e_0=(1-|w|^2)^k e_{w} .
\end{equation}
\begin{equation}\label{r3}
S_k(g)e_{w}=(\delta)^{-2k}e_{g\cdot w}, 
\end{equation}
\end{Remark}

\subsection{The differential action}\label{diff}
We shall suppose that we know the  derived representation $d\pi$ of 
the Lie algebra $\got{g}^J_1$ 
of the Jacobi group $G^J_1$. We
associate to the generators ${{a}}, {{a}}^{\dagger}$ of the Heisenberg group and to the
generators 
$K_{0,+,-}$ of the group  $\text{SU}(1,1)$ the operators ${\mb{a}}, {\mb{a}}^{\dagger}$, respectively 
${\mb{K}}_{0,+,-}$, where 
$({\mb{a}}^{\dagger})^{\dagger}={\mb{a}},~
{\mb{K}}^{\dagger}_0={\mb{K}}_0,
{\mb{K}}^{\dagger}_{\pm}={\mb{K}}_{\mp}$. 
 We take in the definition \eqref{etild} $e_0= e^H_0\otimes
 e^K_0$, i.e. we impose to the cyclic vector $e_0$ to verify simultaneously
 the conditions
\begin{subequations}\label{cond}
\begin{eqnarray}
\mb{a}e_0 & = & 0, \label{a0}\\
~ {\mb{K}}_-e_0 & = & 0, \label{kminus}\\
{\mb{K}}_0e_0 & = & k e_0;~ k>0, 2k=2,3,... .\label{k0}
\end{eqnarray}
\end{subequations}
We consider in  (\ref{k0}) the positive  discrete series
representations $D^+_k$ of $\text{SU}(1,1)$ (cf. \S 9 in \cite{bar47}).

 Perelomov's coherent state vectors   associated to the group $G^J_1$ with 
Lie algebra the Jacobi algebra \eqref{baza},
 based on the Siegel-Jacobi disk  
\begin{equation}\label{mm}
\mc{D}^J_1:= H_1/\R\times
 \text{SU}(1,1)/\text{U}(1)= \C\times\mc{D}_1,
\end{equation}
are defined as 
\begin{equation}\label{csu}
e_{z,w}:=e^{\sqrt{\mu}z{\mb{a}}^{\dagger}+w{\mb{K}}_+}e_0, ~z\in\C,~ |w|<1 .
\end{equation}
The general scheme \cite{sb6} associates to elements of the
Lie algebra $\got{g}$ holomorphic first order  differential operators:
 $X\in\got{g}\rightarrow\db{X}\in\D_1$. The space of
functions on which these operators act in the case of the Jacobi group
 will be made precise later, see \eqref{ofi}.

The following lemma expresses the differential action of the generators
of the Jacobi  algebra as operators of the type $\got{A}_1$,
i.e. first order holomorphic differential operators with polynomial coefficients  in the 
variables $(w,z)$  on $\mc{D}^J_1$, see \cite{jac1}:
\begin{lemma}\label{mixt}The differential action of the generators \eqref{baza}
 of the Jacobi algebra   $\got{g}^J_1$ 
 is given by the formulas:
\begin{subequations}\label{summa}
\begin{eqnarray}
& & \mb{a}=\frac{\pa}{\sqrt{\mu}\pa z};~\mb{a}^{\dagger}=\sqrt{\mu} z+w\frac{\pa}{\sqrt{\mu}\pa z} ;\\
 & & \db{K}_-=\frac{\pa}{\pa w};~\db{K}_0=k+\frac{1}{2}z\frac{\pa}{\pa z}+
w\frac{\pa}{\pa w};\\
& & \db{K}_+=\frac{1}{2}\mu z^2+2kw +zw\frac{\pa}{\pa z}+w^2\frac{\pa}{\pa
w} ,
\end{eqnarray}
\end{subequations}
where $z\in\C$, $|w|<1$.
\end{lemma}

\subsection{The   Jacobi group $G^J_1$}\label{actg}

Following \cite{jac1},  let us introduce the normalized coherent state vectors of the type \eqref{bigch}, based on the
Siegel-Jacobi disk
\begin{equation}\label{psi2}
\underline{e}_{\alpha, w}:=D_{\mu}(\alpha ) S_k(w) e_0;~\alpha\in\C, ~w\in\C, |w|<1. 
\end{equation}

  We find   \cite{jac1} a relation  between the normalized vector
(\ref{psi2}) and the unnormalized Perelomov's CS-vector (\ref{csu}),  
 which play an  important roll  in the proof of Proposition
\ref{mm1}:
\begin{lemma}\label{lema6} The vectors {\em (\ref{psi2}), (\ref{csu})},
i.e.
$$\underline{e}_{\eta, w}:=D_{\mu}(\eta ) S_k(w) e_0;~ e_{z,w'}:=\exp
(\sqrt{\mu}z\mb{a}^{\dagger}+w'{\mb{K}}_+)e_0 . $$
are related by the relation
\begin{subequations}\label{csv}
\begin{align}
\underline{e}_{\eta, w}  & = (1-w\bar{w})^k
\exp (-\frac{1}{2}\mu\bar{\eta}z)e_{z,w}, ~z=\eta-w\bar{\eta},\\ 
e_{z,w} &
=(1-w\bar{w})^{-2k}\exp(\frac{1}{2}\mu\bar{\eta}z)\underline{e}_{\eta,w},
~\eta = \frac{z+\bar{z}w}{1-w\bar{w}}.
\end{align}
\end{subequations}
where the change of variables   $(\eta,w)\mapsto (z,w)$ in \eqref{csv} 
\begin{equation}\label{zwe}
(z,w)= FC(\eta,w), ~ z =
  \eta-w\bar{\eta}
\end{equation} is a $FC$-transform in the
sense of \eqref{etild} for coherent states defined on the
Siegel-Jacobi disk $\mc{D}^J_1$.
\end{lemma}

Let us introduce the notation $\mc{D}^J_1\ni\zeta:=(z,w)\in(\C\times\mc{D}_1)$.
Following the methods of \cite{jac1}, we get  the reproducing kernel 
$K:\mc{D}^J_1\times  \mc{D}^J_1\rightarrow\C$,  $K_{k\mu}(\zeta,\bar{\zeta}')
  :=  (e_{\bar{z},\bar{w}},e_{\bar{z}',\bar{w}'})$:
\begin{equation}\label{ker3}
K_{k\mu}(\zeta,\bar{\zeta}')=
(1-{w}\bar{w}')^{-2k}\exp{\mu F(\zeta,\bar{\zeta}')}, ~
F(\zeta,\bar{\zeta}') = 
\frac{2\bar{z}'{z}+z^2\bar{w}'+\bar{z}'^2w}{2(1-{w}\bar{w}')}.
\end{equation}
In particular, the kernel on diagonal, 
$K_{k\mu}(z,\bar{w})=(e_{\bar{z},\bar{w}},e_{\bar{z},\bar{w}})$, reads
\begin{equation}\label{hot}
K_{ k\mu}(z,\bar{w})=
(1-w\bar{w})^{-2k}\exp{\mu F(z,w)},~ F(z,w)=\frac{2z\bar{z}+z^2\bar{w}+\bar{z}^2w}{2(1-w\bar{w})}, 
\end{equation}
and evidently $K_{ k\mu}(z,\bar{w}) >0,~ \forall (z,w)\in\mc{D}^J_1.$

The holomorphic, transitive and effective  
action of the Jacobi  group
$ G^J_1$ 
on the manifold $\mc{D}^J_1$  (\ref{mm}) can be seen from the
Proposition \ref{mm1}, proved in  \cite{jac1} in the case $\mu =1$:
\begin{Proposition}\label{mm1}
Let us consider the action $S_k(g)D_{\mu}(\alpha )e_{z,w}$, where $g\in
\rm{SU}(1,1)$ has the form \eqref{dg},
$D_{\mu}(\alpha )$ is given by 
 \eqref{deplasare}, and the
coherent state vector $e_{z,w}$ is defined in \eqref{csu}. Then we have the
formula \eqref{xx} and the relations \eqref{xxx}, \eqref{x4}:
\begin{equation}\label{xx}
S_k(g)D_{\alpha}(\alpha )e_{z,w}=\lambda e_{z_1,w_1}, ~ \lambda = \lambda
(g,\alpha; z,w) , \end{equation}
\begin{equation}\label{xxx}
z_1=\frac{\gamma}{\delta}, ~\gamma = z
+\alpha-\bar{\alpha}w;~\delta=\bar{b}w+\bar{a},  ~ w_1=g\cdot
w=\frac{a w+ b}{\delta},\end{equation}
\begin{equation}\label{x4}
\lambda= (\delta)^{-2k}
\exp (-\frac{1}{2}\mu\lambda_1), ~ \lambda_1=\bar{\alpha}(z+\gamma)+\frac{\bar{b}\gamma^2}{\delta}.  
\end{equation}
If we consider the representation
$T_{k\mu}(g,\alpha,t)=S_{k}(g)\beta_{\mu}(\alpha,t)$ of the Jacobi
  group $G^J_1$,
then  $T_{k\mu}(g,\alpha,t)e_{z,w}=\lambda \e^{\ii\mu t} e_{z_1,w_1}$. 

\end{Proposition}

\section{Geometry on the Siegel-Jacobi disk}\label{OFG}
\subsection{The symmetric Fock space}\label{SFS}
In accord with the general scheme of \S \ref{s22}, the scalar product
\eqref{overc} of 
functions from the space $\got{F}_{k\mu}=L^2_{\text{hol}}(\mc{D}^J_1,\rho_{k\mu})$
corresponding to the kernel $K_{k\mu}$ 
defined by \eqref{hot}  on the manifold (\ref{mm}) is:
\begin{equation}\label{ofi}
(\phi ,\psi )_{k\mu}= \int_{\mc{D}^J_1}\!\bar{f}_{\phi}(z,w)f_{\psi}(z,w)\rho_{k\mu}
\dd \nu ,\rho_{k\mu} = \Lambda\! (1\!-\!w\bar{w})^{2k}\!\e^{\!-\mu\frac{2|z|^2+z^2\bar{w}\!+\!\bar{z}^2w}{2(1\!-\!w\bar{w})}}.
\end{equation}
where the value of the $G^J_1$-invariant measure $\dd\nu $ 
\begin{equation}\label{ofi3}
d \nu =\mu\frac{\dd \Re w \dd\Im w}{(1-w\bar{w})^3}\dd \Re z \dd \Im z
\end{equation}
is  obtained  in  (\ref{dnu}).  
Note that the value of the normalization constant in \eqref{ofi} is the same as in the case $\mu=1$, i.e.
\begin{equation}\label{ofi1}
\Lambda = \frac{4k-3}{2\pi^2} .
\end{equation}

The base of orthonormal functions associated to the CS-vectors attached  to the
Jacobi group $G^J_1$, defined  on the manifold $\mc{D}^J_1$ (\ref{mm})
\cite{jac1}, \cite{sbcg}
\begin{equation}\label{x3x1}
f_{{|n>;e_{k',k'+m}}}(z,w)  := 
(e_{\bar{z},\bar{w}},|n> e_{k',k'+m}),~z\in\C,~ |w|<1,
k=k'+\frac{1}{4},2k'\in \Z_+
\end{equation}
consists 
of the holomorphic polynomials 
\begin{equation}\label{x4x}
f_{{|n>;e_{k',k'+m}}}(z,w)=f_{e_{k',k'+m}}(w)\frac{P_n(\sqrt{\mu}z,w)}{\sqrt{n!}}, 
\end{equation}
where the monomials  $f_{e_{k',k'+m}}$ are defined in \eqref{f2},
while   
\begin{equation}\label{marea}
P_n(z,w)=n!\sum _{p=0}^{[\frac{n}{2}]}
(\frac{w}{2})^p\frac{z^{n-2p}}{p!(n-2p)!} .
\end{equation}
 
The series expansion \eqref{sumBERG}  of the  reproducing kernel \eqref{ker3} 
 reads
\begin{equation}
\label{ker31}K_{k\mu}(z,w;\bar{z}',\bar{w}')   =
\sum_{n,m}f_{|n>,e_{k',k'+m}}(z,w)\bar{f}_{|n>,e_{k',k'+m}}(z',w').
\end{equation}
\subsection{Two-forms}
We follow further  the general prescription of \S\ref{s22}.
We calculate the K\"ahler potential on $\mc{D}^J_1$ 
as the logarithm of the reproducing kernel 
$f:=\ln K_{k\mu}$, 
i.e.
\begin{equation}\label{keler}
f =\mu\frac{2z\bar{z}+z^2\bar{w}+\bar{z}^2w}{2(1-w\bar{w})}
-2k\ln (1-w\bar{w}) .
\end{equation}
The K\"ahler two-form $\omega$ is obtained  with 
 formulas \eqref{kall}, \eqref{Ter}, i.e. 
\begin{equation}\label{aaa}
-\ii \omega = h_{z\bar{z}}dz\wedge d\bar{z}+h_{z\bar{w}}dz\wedge
d\bar{w}
-h_{\bar{z}w}d\bar{z}\wedge d w +h_{w\bar{w}}dw\wedge d\bar{w} .
\end{equation}
The volume form is:
\begin{equation}\label{vol1}-\omega\wedge\omega=2\left|\begin{array}{cc}
h_{z\bar{z}} & h_{z\bar{w}}\\
h_{\bar{z}w} & h_{\bar{w}w}
\end{array}\right|
\dd z\wedge \dd \bar{z}\wedge \dd w \wedge \dd \bar{w} .
\end{equation}
With  the K\"ahler potential \eqref{keler}, we get the metric
coefficients  in \eqref{aaa}
\begin{equation}\label{metrica}
h_{z\bar{z}}=\frac{\mu}{P}, ~h_{z\bar{w}}=\mu\frac{\eta}{P},~
  h_{w\bar{w}}=\mu\frac{|\eta|^2}{P}+\frac{2k}{P^2}, ~P=(1-w\bar{w}),
\end{equation}
where $\eta$ was defined previously in \eqref{csv}.

If  we introduce the notation 
\begin{equation}\label{DETG}
G(z):=\det
(h_{\alpha\bar{\beta}})_{\alpha,\beta=1,\dots,n},\end{equation}
 then
the {\it  Ricci form of the Bergman metric} is  (see p. 90 in  \cite{mor})
\begin{equation}\label{RICCI}
\rho_M(z):=\ii \sum_{\alpha,\beta=1}^n\text{Ric}_{\alpha\bar{\beta}}(z)\dd
z_{\alpha}\wedge\dd \bar{z}_{\beta}, ~ \text{Ric}_{\alpha\bar{\beta}}(z)=
 -\frac{\pa^2}{\pa z_{\alpha}\pa \bar{z} _{\beta}} \ln G(z).
\end{equation}

For the determinant $G$ \eqref{DETG} for the Siegel-Jacobi disk we
find \begin{equation}G(z,w)=\frac{2k\mu}{(1-w\bar{w})^3}.\end{equation} 

The scalar curvature at a point $p\in M$ of coordinates $z$ is (see
p. 294 in \cite{koba1})
\begin{equation}\label{scc}
s_M(p)=\sum_{\alpha,{\beta}=1}^n
(h_{\alpha{\bar\beta}})^{-1}\text{Ric}_{\alpha\bar{\beta}}(z). 
\end{equation}

Following \cite{LU08},  let us introduce also the positive definite
(1,1)-form on $M=G/H$
\begin{equation}\label{RIC2}\tilde{\omega}_M(z):= 
\ii\sum_{\alpha,\beta\in\Delta_+} \tilde{h}_{\alpha\bar{\beta}} (z) \dd z_{\alpha}\wedge
\dd\bar{z}_{\beta}, ~  \tilde{h}_{\alpha\bar{\beta}} (z):= 
  (n+1)h_{\alpha\bar{\beta}} (z)-  \text{Ric}_{\alpha\beta}(z),
\end{equation}
which is K\"ahler, with K\"ahler potential $\tilde{f} = \ln
(K(z)^{n+1}G(z))$.

With \eqref{keler}-\eqref{RIC2}, we obtain
\begin{Proposition}\label{prop2}
The K\"ahler two-form $\omega_{k\mu}$ on $\mc{D}^J_1$, $G^J_1$-invariant to the
action \eqref{xxx} is:
\begin{equation}\label{aab}
-\ii\omega_{k\mu}(z,w) =2k\frac{\dd w \wedge \dd\bar{w}}{P^2} +
\mu\frac{A\wedge \bar{A}}{P},
~A=\dd z+\bar{\eta}\dd w, ~\eta=\frac{z+\bar{z}w}{P}.
\end{equation}
The  hermitian metric on $G^J_1$ corresponding to the K\"ahler two-form
\eqref{aab}
is \begin{equation}\label{metric}
\dd s^2_{k\mu}(z,w)=2k\frac{\dd w \otimes \dd\bar{w}}{P^2} +
\mu\frac{A\otimes\bar{A}}{P}.
\end{equation}

The volume form \eqref{vol1} is:
\begin{equation}\label{dnu}
\omega_{k\mu}\wedge\omega_{k\mu} = 4k\mu (P)^{-3} 4\dd \Re z \wedge \dd \Im z \wedge\dd
 \Re w  \wedge \dd \Im w  .
\end{equation}

The K\"ahler two-form \eqref{RIC2} for $\mc{D}^J_1$, corresponding to
the K\"ahler potential $$\tilde{f}(z,w)= 3[\mu f(z,w)-(2k+1)\ln
(1-w\bar{w})] $$ reads
\begin{equation}\label{doiD}
-\ii \tilde{\omega}_{\mc{D}^J_1} = 
3 \left[(2k+1)\frac{\dd w \wedge \dd\bar{w}}{P^2} +
\mu\frac{A\wedge \bar{A}}{P}\right].
\end{equation}
The Ricci form \eqref{RICCI} is
\begin{equation}\label{RICJD}
\rho_{\mc{D}^J_1}(z,w) = -3\ii\frac{\dd w\wedge\dd \bar{w}}{(1-w\bar{w})^2},
\end{equation}
and  $\mc{D}^J_1$   is not an Einstein manifold with respect to the
metric \eqref{metric}.

 The scalar curvature \eqref{scc}   has the value
\begin{equation}\label{sCAL}
s_{\mc{D}^J_1}(p)=-\frac{3}{2k}, p\in \mc{D}^J_1.\end{equation}
\end{Proposition}
We emphasize that  result \eqref{sCAL} was previously obtained in \cite{jae}.

Now we introduce  in the expression of the one-form  $A$ in
\eqref{aab}  the $FC$-transform \eqref{zwe}, $z=\eta-w\bar{\eta}$, so $A= \dd
\eta - w\dd \bar{\eta}$, and then $A\wedge\bar{A}= P\dd \eta\wedge
\dd\bar{\eta}$.
We  underline   the significance of the change of coordinates $FC$   which realizes 
{\it the  fundamental conjecture for the
Siegel-Jacobi disk},  proved in Proposition 3.1  in \cite{FC}:
\begin{Proposition}\label{prop3}
The  $FC$-transform \eqref{zwe},  $FC(\eta,w)=(z,w)$, which appears in
{\emph{Le\-mma \ref{lema6}}},
$z=\eta-w\bar{\eta}$, is a  homogeneous K\"ahler diffeomorphism,
i.e.  $FC^*\omega_{k\mu}(z,w)= \omega_{k\mu}(\eta,w)$, where 
\begin{equation}\label{aabc}
\omega_{k\mu}(\eta,w) =\omega_{k}(w)+\omega_{\mu}(\eta).
\end{equation}
The   K\"ahler two-form \eqref{aabc}  is invariant  to the action of
$G^J_1$ on $\C\times\mc{D}_1$, $((g,\alpha), (\eta,w))\rightarrow
(\eta_1,w_1)$, where
\begin{equation}\label{grozav}
\eta_1=a(\eta+\alpha)+b(\bar{\eta}+\bar{\alpha}), ~
w_1=\frac{aw+b}{\bar{b}w+\bar{a}}, ~g=\left(\begin{array}{cc} a & b\\
      \bar{b}& \bar{a}\end{array}\right)\in {\emph{\text{SU}}}(1,1). 
\end{equation}
\end{Proposition}

\subsection{Geodesics}\label{GEOS}
The equations of geodesics on a manifold $M$ with linear connection
with components of the linear connections $\Gamma$ are  (see
e.g. Proposition 7.8 p. 145 in \cite{koba1}) 
\begin{equation} \label{GEOO}
\frac{\dd ^2 x_i}{\dd t^2}  +\sum_{j,k}\Gamma^i_{jk}\frac{\dd x_j}{\dd t}
\frac{\dd x_k}{\dd t}  =  0, ~i=1,\dots,n .
\end{equation}
 $\dd s^2(z,\bar{z})$ from \eqref{herm} gives  the hermitian metric of $M$ in
local coordinates, while the fact the metric is K\"ahlerian imposes
the restrictions \eqref{condH}.
The non-zero Christoffel's  symbols  $\Gamma$ are determined by the equations
(cf. (12) at p. 156  in \cite{koba2})  
\begin{equation}\label{CRISTU} \sum_{\alpha} h_{\alpha\bar{\epsilon}}\Gamma^{\alpha}_{\beta\gamma}=
\frac{\pa  h_{\bar{\epsilon}\beta}}{\pa z_{\gamma}}. 
\end{equation}
In the variables $(z,w)\in(\C,\mc{D}_1)$ the equations of geodesics
\eqref{GEOO}  read
\begin{equation}
 \left\{
 \begin{array}{l}
 \frac{\dd^2 z}{\dd t^2}+\Gamma^z_{zz}\left(\frac{\dd z}{\dd
     t}\right)^2 +
2\Gamma^z_{zw}\frac{\dd z}{\dd t}  \frac{\dd w}{\dd t} +\Gamma
^z_{ww}\left(\frac{\dd w}{\dd t}\right)^2 =0   ;\\
  \frac{\dd^2 w}{\dd t^2}+\Gamma^w_{zz}\left(\frac{\dd z}{\dd
     t}\right)^2 +
2\Gamma^w_{zw}\frac{\dd z}{\dd t}  \frac{\dd w}{\dd t} +\Gamma
^w_{ww}\left(\frac{\dd w}{\dd t}\right)^2=0 .
     \end{array}
 \right.
\end{equation}
The equations \eqref{CRISTU} which determine the $\Gamma$-symbols for
the Siegel-Jacobi disk are 
\begin{equation}\label{ec528}
 \left\{
 \begin{array}{l}
h_{z\bar{z}}\Gamma^z_{zz}+h_{w\bar{z}}\Gamma^w_{zz}=
\frac{\pa h_{z\bar{z}}}{\pa z}; \\
h_{z\bar{w}}\Gamma^z_{zz}+h_{w\bar{w}}\Gamma^w_{zz}=
\frac{\pa h_{z\bar{w}}}{\pa z}. 
\end{array}
 \right.
\end{equation}
\begin{equation}\label{ec529}
 \left\{
 \begin{array}{l}
h_{z\bar{z}}\Gamma^z_{zw}+h_{w\bar{z}}\Gamma^w_{zw}=
\frac{\pa h_{w \bar{z}}}{\pa z}; \\
h_{z\bar{w}}\Gamma^z_{zw}+h_{w\bar{w}}\Gamma^w_{zw}=
\frac{\pa h_{w\bar{w}}}{\pa z}. 
\end{array}
 \right.
\end{equation}
\begin{equation}\label{ec530}
\left\{
 \begin{array}{l}
h_{z\bar{z}}\Gamma^z_{ww}+h_{w\bar{z}}\Gamma^w_{ww}=
\frac{\pa h_{w \bar{z}}}{\pa w}; \\
h_{z\bar{w}}\Gamma^z_{ww}+h_{w\bar{w}}\Gamma^w_{ww}=
\frac{\pa h_{w\bar{w}}}{\pa w}. 
\end{array}
 \right.
\end{equation}
With \eqref{metrica}, we calculate  easily the derivatives
\begin{equation}\label{ec531}
\begin{split}
\frac{\pa h_{z\bar{z}}}{\pa z } &= 0;~ \frac{\pa h_{z\bar{w}}}{\pa z}=
\frac{\mu}{P^2}: ~\frac{\pa h_{w\bar{z}}}{\pa z}=
\mu\frac{\bar{w}}{P^2}; \frac{\pa h_{w\bar{w}}}{\pa z }=\mu\frac{\bar{\eta}+\eta\bar{w}}{P^2};\\
\frac{\pa h_{w\bar{z}}}{\pa w } & = 2\mu\frac{\bar{w}\bar{\eta}}{P^2}; \frac{\pa h_{w\bar{w}}}{\pa w}=
\mu\frac{\bar{z}\bar{\eta}}{P^2}+3\mu\frac{\bar{w}|\eta|^2}{P^2}+4k\frac{\bar{w}}{P^3}.
\end{split}
\end{equation}
Introducing \eqref{ec531} into \eqref{ec528}-\eqref{ec530}, we find
for   the  Christoffel's symbols $\Gamma$-s  the expressions
\begin{equation}
\begin{split}
\Gamma^z_{zz}  & =-\lambda\bar{\eta},~\Gamma^w_{zz}=\lambda; ~
\Gamma^z_{zw}=-\lambda\bar{\eta}^2+\frac{\bar{w}}{P};\\
 ~\Gamma^w_{wz} & =\lambda\bar{\eta};~
 \Gamma^z_{ww}=-\lambda\bar{\eta}^3; ~ \Gamma^w_{ww} =
 \lambda\bar{\eta}^2+2\frac{\bar{w}}{P}, ~ \lambda = \frac{\mu}{2k},
\end{split}
\end{equation}
and we can formulate
\begin{Remark}\label{GER}
The equations of geodesics on the Siegel-Jacobi corresponding to
metric defined by $\omega_{k\mu}$  \eqref{aab} are
\begin{subequations}\label{geo}
\begin{align}
\mu\bar{\eta}G^2_1 & = 2k G_3, ~ G_1=\frac{\dd z}{\dd
  t}+\bar{\eta}\frac{\dd w}{\dd t}, ~ G_3 = \frac{\dd ^2z}{\dd
  t^2}+2\frac{\bar{w}}{P}\frac{\dd z}{\dd t}\frac{\dd w}{\dd t}, ~ P=1-w\bar{w};  \label{geo1}\\
\mu G^2_1 & = -2k G_2, ~G_2=\frac{\dd^2w}{\dd
  t^2}+2\frac{\bar{w}}{P}(\frac{\dd w}{\dd t})^2 . \label{geo2}
\end{align}
\end{subequations}
\end{Remark}
The above equations were written in the case
$\mu=1$ in \cite{jac1}.

Now we discuss the solution of the system \eqref{geo}.

a) If in the system of differential equations \eqref{geo}, we take
$\mu=0$, then   from \eqref{geo2} we get $G_2=0$, i.e.  the equation of
geodesics on $\mc{D}_1$. The solution of the equation
$G_2=0$, with $
w(0)=0$,  $\dot{w}(0)=B$, is given by the $FC$ transform \eqref{u5}, i.e. 
\begin{equation}\label{solg}
w(t)=\frac{B}{|B|}\tanh(t|B|) . 
\end{equation} 
If we introduce the solution \eqref{solg} of the equation $G_2=0$ of  geodesics on $\mc{D}_1$,
 then  the solution  of the differential equation $G_3=0$ with the
initial condition $z_0=\frac{\dd z}{\dd t}|_{t=0}$, $z_1= z(0)$ is
$z(t)=\frac{z_0}{B}w(t)+z_1$, and \eqref{geo1} is satisfied.

b) If $\mu\not= 0$,  a {\it particular solution  $(z,w)$  of
  \eqref{geo} of the system of geodesics on the
Siegel-Jacobi disk $\mc{D}^J_1$ is given by 
$z=\eta_0-\bar{\eta}_0w(t)$,   and $w=w(t)$   with $w(t)$ given by \eqref{solg}, and $
\eta=\eta_0$ independent of $t$}. This particular solution has  been noted already in \cite{FC}   in the case $\mu=1$.  This is
{\it a particular case of the solution $\eta=\eta_0 + t\eta_1$ of equation
of geodesics $\frac{\dd ^2 \eta}{\dd t^2}=0$
 on the flat manifold $\C$}
corresponding to the separation of variables as in \eqref{aabc} of
$\omega_{k\mu}(\eta,w)$.  We can formulate this observation:
\begin{Remark}\label{part}
The $FC$ transform given by \eqref{zwe} in {\emph {Lemma  \ref{lema6} }}is not a
geodesic mapping, but it gives
geodesics $(z(t),w)=FC(\eta,w)$ on the
non-symmetric space $\mc{D}^J_1$ with $w=w(t)$, given
by \eqref{solg}, and $ \eta =\eta_0$.
\end{Remark}

\subsection{Embeddings}\label{EMBSS}

We recall that the homogeneous K\"ahler manifolds  $M=G/H$  which admit an
embedding of the type  given by Remark \ref{HTR} are called
{\it coherent type manifolds}, and the groups $G$ are called
{coherent-state type groups} \cite{lis,neeb}. We particularize Remark
\ref{HTR} in the case of the Siegel-Jacobi disk and we have 
\begin{Remark} \label{REM9}
The Jacobi group $G^J_1$
is a coherent-state type group and the Siegel-Jacobi disk  $\mc{D}^J_1$ is a quantizable K\"{a}hler
 coherent
state manifold. The Hilbert space of functions 
\fl ~ is the space $ \got{F}_{k\mu}=
L^2_{\text{hol}}(\mc{D}^J_1,\rho_{k\mu})$ with  the scalar product
\eqref{ofi}-\eqref{ofi3}. 
The   K\"ahlerian embedding  $\iota_{\mc{D}^J_1}:\mc{D}^J_1\hookrightarrow
\db{CP}^{\infty}$  \eqref{invers}  $\iota_{\mc{D}^J_1} =[\Phi]=
[\varphi_0:\varphi_1:\dots\varphi_N: \dots ] $ is realized with  an ordered version of the base functions
$\Phi =\left\{ f_{{|n>;e_{k',k'+m}}}(z,w) \right\}$  given by \eqref{x4x}, 
 and the K\"ahler two-form  \eqref{aab} is the pullback of
the Fubini-Study K\"ahler two-form \eqref{FBST} on $\db{CP}^{\infty}$,
$$\omega_{k\mu}= \iota_{\mc{D}^J_1}^*\omega_{FS}|_{\db{CP}^{\infty}},
~\omega_{k\mu}(z,w)=\omega_{FS}([\varphi_N(z,w)]).$$ 
The normalized Bergman kernel \eqref{kmic} of the Siegel-Jacobi disk $\kappa_{k\mu}$
expressed in the variables $\zeta=(z,w)$, $\zeta'=(z',w')$ reads
\begin{equation}\label{redBKJ}
\kappa_{k\mu}(\zeta,\bar{\zeta}')=\kappa_{k}(w,\bar{w}')\exp[\mu
(F(\zeta,\bar{\zeta}')-\frac{1}{2}(F(\zeta)+F(\zeta'))], 
\end{equation}
where $\kappa_{k}(w,\bar{w}')$ is the normalized Bergman kernel for
the Siegel disk $\mc{D}_1$
\begin{equation}
\kappa_{k}(w,\bar{w}')=\left[\frac{(1-|w|^2)(1-|w'|^2)}{(1-w\bar{w}')^2}\right]^k,
\end{equation}
 $F(\zeta,\bar{\zeta}')$ is defined in \eqref{ker3}, and
$F(\zeta)$ is defined in \eqref{hot}. 
The Berezin  kernel of $\mc{D}^J_1$ is 
$$b_{k\mu}(\zeta,\zeta')=b_{k}(w,w')\exp\mu[2\Re
F(\zeta,\bar{\zeta}')-F(\zeta)- F(\zeta')],$$
where $b_{k}(w,w')=|\kappa_{k}(w,\bar{w}')|^2$.

With formula \eqref{ker31}, we get for the
diastasis function on the Siegel-Jacobi disk the expression:
\begin{equation}\label{dia1}
\frac{D_{k\mu}(\zeta,\zeta')}{2}=k
  \ln\frac{\vert 1-w\bar{w}'\vert^2}{(1-|w|^2)(1-|w'|^2)} +\mu[\frac{F(\zeta)+F(\zeta')}{2}-\Re F(\zeta,\bar{\zeta'})].
\end{equation}
\end{Remark}
\section{Appendix}
We recall some definitions from differential geometry of notions used
in the paper.

\begin{deff}({\bf Reductive homogeneous spaces}, cf. Nomizu \cite{nomizu})
A homogeneous space $M=G/H$ is {\em reductive} if the Lie algebra \g~
of $G$ may be decomposed into a vector space direct sum of the Lie
algebra \h ~ of $H$ and an $\Ad (H)$-invariant subspace \m , that is
\begin{subequations} \label{dec}
\begin{equation}\label{sum1}
\g = \h + \m, ~ \h\cap\m =0,
\end{equation}
\begin{equation}\label{sum2}
 \Ad(H)\m \subset \m. 
\end{equation}
Condition \eqref{sum2} implies
\begin{equation}\label{sum3}
[\h ,\m ]\subset \m
\end{equation}
\end{subequations}
and, conversely, if $H$ is connected, then \eqref{sum3} implies
 \eqref{sum2}. Note that $H$ is always connected if $M$ is simply
connected. The decomposition \eqref{sum1} verifying \eqref{sum2} is called
a  $H$-stable decomposition.
\end{deff}

We recall that an analytic  Riemannian manifold is called {\it
  Riemannian globally symmetric}  if each $p\in M$ is an isolated
fixed point of an involutive isometry $s_p$ of $M$, see p. 205 in \cite{helg}.
We recall that if the Lie algebra  $\got{g}$ and its  subalgebra $\got{h}$
associated with  the homogeneous manifold $M=G/H$  satisfy \eqref{dec},
then a necessary and sufficient condition for $M$  to be a  {\it
 locally  symmetric space}  is \begin{equation}\label{symsp}
[\got{h},\got{h}]\subset\got{h},~ [\got{h},\got{m}]\subset\got{m},~
[\got{m},\got{m}]\subset\got{h}. 
\end{equation}
{\it If $M$ is a complete, simply connected Riemannian locally symmetric
  space, then $M$ is a Riemannian globally symmetric space} (cf. Theorem 5.6
p 232 in \cite{helg}).

Let  $\text{Exp}_p:M_p\rightarrow M$ be  the geodesic exponential  map
from the tangent space  $M_p$ at $p\in M$ to $M$ (cf. definition at
p. 33 in \cite{helg}), and also let $\exp:\got{g}\rightarrow G$ be  the
exponential map from the Lie algebra $\got{g}$ to the Lie group $G$.  We also recall  that the symmetric spaces have the property that 
 $\text{Exp}_p(\got{m})$$=\exp_p(\got{m})$, i.e. the exponential $\exp_p$ from
 the Lie algebra to the Lie group at $p$ gives geodesics such that
 lines in the tangent space are geodesics  in the manifold
 emerging from the point $p$. More
exactly, if $\sigma$ is the projection $\sigma:G\rightarrow G/H$,
such that $o=\sigma(e)$, where $e$ is the unit element in the group
$G$, then for symmetric spaces $\sigma(\exp X)=\text{Exp}(\dd \sigma
X), X\in\got{m}$, cf. Theorem 3.3 p 212 in 
\cite{helg}.

We also have been used in Remark \ref{part} the notion of {geodesic
  mapping} (cf. Definition 5.1 p. 127 in \cite{mikes})
\begin{deff}
If $f:M\rightarrow N$ is a diffeomorphism between manifolds with
affine connections, then $f$ is called a geodesic mapping if it maps
geodesic curves on $M$ into geodesic curves on $N$.
\end{deff}
\enlargethispage{1cm}
\subsection*{Acknowledgments}
This research  was  supported by the 
ANCS project  program PN 09 37 01 02/2009 and  by the UEFISCDI - Romania
 program PN-II Contract No. 55/05.10.2011. I am grateful to Mircea Bundaru for the suggestions to improve the manuscript. The author  would like to express his thanks to Professors  Jae-Suk Park, 
Jae-Hyun Yang  and  Kisik Kim for the possibility to present a
preliminary version of this paper in seminaries at 
the Center for Geometry and Physics, 
Institute for Basic Science, Pohang, Korea, Department of Mathematics
and, respectively,  Department of Physics,  Inha University
Incheon, Korea. The author is
 indebted to the unknown refrees for the sugesstions and criticism.

\today
\end{document}